\documentclass[journal,twoside,web]{ieeecolor}
\usepackage{lcsys}
\usepackage{cite}
\usepackage{amsmath,amssymb,amsfonts}
\usepackage{algorithmic}
\usepackage{graphicx}
\usepackage{textcomp}
\usepackage{multirow}
\usepackage{bm}

\newcommand{\ignore}[1]{}

\newcommand{\norm}[1]{\left\Vert#1\right\Vert} 

\newcommand{\bbm}{\begin{bmatrix}}
\newcommand{\ebm}{\end{bmatrix}}
\newcommand{\bma}[1]{\left[\begin{array}{#1}}
\newcommand{\ema}{\end{array}\right]}

\DeclareMathAlphabet{\mbf}{OT1}{ptm}{b}{n}
\newcommand{\mbs}[1]{{\boldsymbol{#1}}}
\newcommand{\mbc}[1]{ \boldsymbol{\mathcal{#1}} } 

\newcommand{\mbsdot}[1]{{\dot{\boldsymbol{#1}}}}

\newcommand{\mbfdot}[1]{{\dot{\mbf{#1}}}}
\newcommand{\mbfbar}[1]{{\bar{\mbf{#1}}}}
\newcommand{\mbfhat}[1]{{\hat{\mbf{#1}}}}

\newcommand{\mbftilde}[1]{{\tilde{\mbf{#1}}}}

\def\fdotb{{\raisebox{-0.6ex}{ \kern0.2ex\raisebox{0.8ex}{\tiny $\hspace*{-1ex}\circ$}}}}

\def\fddotb{{\raisebox{-0.6ex}{ \kern0.2ex\raisebox{0.8ex}{\tiny $\hspace*{-1ex}\circ\circ$}}}}

\newcommand{\f}{\frac}

\newcommand{\trans}{{\ensuremath{\mathsf{T}}}} 
\newcommand{\beq}{\begin{equation}}
\newcommand{\eeq}{\end{equation}}
\newcommand{\bdis}{\begin{displaymath}}
\newcommand{\edis}{\end{displaymath}}
\newcommand{\beqarray}{\begin{eqnarray}}
\newcommand{\eeqarray}{\end{eqnarray}}
\newcommand{\beqarraynn}{\begin{eqnarray*}}
\newcommand{\eeqarraynn}{\end{eqnarray*}}
\newcommand{\eye}{\mbf{1}}
\newcommand{\inv}{^{-1}}

\def\BibTeX{{\rm B\kern-.05em{\sc i\kern-.025em b}\kern-.08em
		T\kern-.1667em\lower.7ex\hbox{E}\kern-.125emX}}
\begin{document}
%
%
%
%
%
%
%
\def \myJournal {IEEE Control Systems Letters}
\def \myDoi {10.1109/LCSYS.2023.3237954}
\def \myPaperSiteName {IEEE Xplore}
\def \myPaperSiteLink {https://ieeexplore.ieee.org/document/10020153}
\def \myYear {2023}
\def \myPaperCitation{F. Ahmed, L. Sobiesiak and J. R. Forbes, ``Cascaded Model Predictive Control of a Tandem-Rotor Helicopter,'' in \textit{IEEE Control Systems Letters}, vol. 7, pp. 1345-1350, 2023.}


\begin{figure*}[t]

\thispagestyle{empty}
\begin{center}
\begin{minipage}{6in}
\centering
This paper has been accepted for publication in \emph{\myJournal}. 
\vspace{1em}

This is the author's version of an article that has, or will be, published in this journal or conference. Changes were, or will be, made to this version by the publisher prior to publication.
\vspace{2em}

\begin{tabular}{rl}
DOI: & \myDoi\\
\myPaperSiteName: & \texttt{\myPaperSiteLink}
\end{tabular}

\vspace{2em}
Please cite this paper as:

\myPaperCitation

\vspace{15cm}
\copyright \myYear \hspace{4pt}IEEE. Personal use of this material is permitted. Permission from IEEE must be obtained for all other uses, in any current or future media, including reprinting/republishing this material for advertising or promotional purposes, creating new collective works, for resale or redistribution to servers or lists, or reuse of any copyrighted component of this work in other works.

\end{minipage}
\end{center}
\end{figure*}
\newpage
\clearpage
\pagenumbering{arabic} 
	
\title{Cascaded Model Predictive Control of a Tandem-Rotor Helicopter}
\author{Faraaz Ahmed, Ludwik Sobiesiak, and James Richard Forbes
\thanks{This work is supported by NGC Aerospace as well as the Mitacs Accelerate and NSERC Discovery Grant programs.}
\thanks{Faraaz Ahmed (faraaz.ahmed@mail.mcgill.ca) and James Richard Forbes (james.richard.forbes@mcgill.ca) are with the Department of Engineering, McGill University, Montreal QC, Canada, H3A 0C3.}
\thanks{Ludwik Sobiesiak (ludwik.sobiesiak@ngcaerospace.com) is with the Engineering Department at NGC Aerospace Ltd., Sherbrooke QC, Canada, J1L 2T9.}}

\maketitle
\thispagestyle{empty}

\begin{abstract}
	This letter considers cascaded model predictive control (MPC) as a computationally lightweight method for controlling a tandem-rotor helicopter. A traditional single MPC structure is split into separate outer and inner-loops. The outer-loop MPC uses an $SE_2(3)$ error to linearize the translational dynamics about a reference trajectory. The inner-loop MPC uses the optimal angular velocity sequence of the outer-loop MPC to linearize the rotational dynamics. The outer-loop MPC is run at a slower rate than the inner-loop allowing for longer prediction time and improved performance. Monte-Carlo simulations demonstrate robustness to model uncertainty and environmental disturbances. The proposed control structure is benchmarked against a single MPC algorithm where it shows significant improvements in position and velocity tracking while using significantly less computational resources.
\end{abstract}


\vspace{-0.5em}
\section{Introduction}
\IEEEPARstart{T}{he} tandem-rotor helicopter is an unconventional unmanned aerial vehicle (UAV) design that has two main-rotors mounted along the longitudinal axis. This design offers several advantages, including a large center of mass range, and larger lift capacities with smaller rotors \cite{Johnson2013}.

Several control methods have been explored in the literature. In \cite{Downing1987}, a proportional-integral controller is developed for a Chinook autoland system. A nonlinear Lyapunov algorithm for hover control is demonstrated in \cite{Dzul2002}, and an $\mathcal{L}_1$ adaptive controller is designed using a linearized dynamics in \cite{Gaoyuan2018}. Model predictive control (MPC) of a tandem-rotor helicopter using linearized system dynamics is presented in \cite{Ahmed2022}. Although these methods are effective in simulation, smaller UAV platforms are often limited in computational power making real-time implementation impractical.

The objective of this letter is to present a computationally lightweight control algorithm that is better suited for real-time implementation. MPC is selected for its ability to easily handle constraints and the wide availability of efficient quadratic program (QP) solvers. The single MPC (SMPC) strategy of \cite{Ahmed2022} results in a large control structure. The size and complexity of the QP used in the MPC solution grows quickly with the state dimension and horizon length. Combining a large QP and fast update rate is especially problematic.

In this letter, a cascaded MPC (CMPC) approach is considered where the single controller from \cite{Ahmed2022} is split into separate independent outer-loop and inner-loop MPC algorithms. CMPC has been used to effectively control multi-timescale systems, such as in \cite{Schlagenhauf2020} where cascaded nonlinear MPC (NMPC) is demonstrated on a quadrotor, and in \cite{Ulbig2011} where a linear CMPC structure is applied to a power system. Here, the CMPC strategy provides a method of effectively handling the multi-timescale nonlinear dynamics of the tandem-rotor helicopter in a less computationally demanding way than the SMPC structure of \cite{Ahmed2022} and the NMPC structure of \cite{Schlagenhauf2020}. Unlike \cite {Ulbig2011}, which deals with linear time-invariant (LTI) dynamics, this approach is applied to linearized dynamics that are linear time-varying (LTV). Although this necessitates solving two QPs in real-time, each MPC has fewer states and optimization variables, and can be solved faster.

The novel contribution of this letter is the synthesis of a multi-rate CMPC structure, for nonlinear UAV control, where the outer-loop dynamics are linearized about a reference trajectory, and the inner-loop dynamics are continuously linearized about the outputs of the outer-loop controller. The resulting formulation is less computationally complex than the SMPC structure of \cite{Ahmed2022} and the NMPC structure of \cite{Schlagenhauf2020}, while still providing good tracking performance.

\vspace{-0.5em}
\section{Kinematics and Dynamics}
The standard North-East-Down basis vectors are used to define an inertial frame $\mathcal{F}_a$ \cite{Hughes2004}. Let point $w$ represent a point in $\mathcal{F}_a$. Denote $\mathcal{F}_b$ as a frame fixed to the body of the helicopter. The direction cosine matrix (DCM) $\mbf{C}_{ab} \in SO(3)$ relates the attitude of $\mathcal{F}_a$ to $\mathcal{F}_b$. A physical vector $\underrightarrow{v}$ can be resolved in either $\mathcal{F}_a$ as $\mbf{v}_a$, or in $\mathcal{F}_b$ as $\mbf{v}_b$, where $\mbf{v}_a = \mbf{C}_{ab} \mbf{v}_b$.

Consider a tandem-rotor helicopter \cite{Dzul2002}, with mass $m_\mathcal{B}$, modeled as a rigid body with thrust, gravitational, and drag forces acting on it. Let $z$ be a point collocated with the center of the mass of the helicopter. The translational and rotational kinematics are given by
\begin{subequations}\label{eq201}
	\begin{align}
		\mbfdot{C}_{ab} &= \mbf{C}_{ab}\mbs{\omega}^{ba^\times}_b, \\
		\mbfdot{r}^{zw}_a &= \mbf{v}^{zw/a}_a, 
	\end{align}
\end{subequations}
where $\mbf{r}^{zw}_a$ is the position of point $z$ relative to point $w$, resolved in $\mathcal{F}_a$, $\mbf{v}^{zw/a}_a$ is the velocity of point $z$ relative to point $w$ with respect to $\mathcal{F}_a$, resolved in $\mathcal{F}_a$, and $\omega^{ba}_b$ is the angular velocity of $\mathcal{F}_b$ relative to $\mathcal{F}_a$. The cross operator $(\cdot)^\times$ is a mapping from $\mathbb{R}^3$ to the Lie algebra $\mathfrak{so}(3)$ such that $\mbf{u}^\times \mbf{v} = -\mbf{v}^\times \mbf{u}$.

The translational and rotational dynamics are \cite{Faessler2018}
\begin{subequations}\label{eq202}
	\begin{align}
		\hspace{-10pt}
		m_\mathcal{B}\mbfdot{v}^{zw/a}_a &= \eye_3 m_\mathcal{B} g - \mbf{C}_{ab} \eye_3 f - \mbf{C}_{ab}\mbf{D}\mbf{C}_{ab}^\trans\mbf{v}^{zw/a}_a, \\
		\mbf{J}^{\mathcal{B}z}_b \mbsdot{\omega}^{ba}_b &= \mbf{m}_b - \mbf{E}\mbf{C}_{ab}^\trans\mbf{v}^{zw/a}_a - \mbf{F}\mbs{\omega}^{ba}_b - \mbs{\omega}^{ba^\times}_b \mbf{J}^{\mathcal{B}z}_b \mbs{\omega}^{ba}_b,
	\end{align}
\end{subequations}
where $\mbf{J}^{\mathcal{B}z}_b$ is the helicopter's second moment of mass resolved in $\mathcal{F}_b$, $f$ is the total thrust force from the rotors, $\mbf{m}_b$ is the total control torque from the rotors, $g = 9.81 \text{m/s}^2$, and $\eye_3 = [0\;0\;1]^\trans$. The matrix $\mbf{D} = \text{diag}(d_x,\;d_y\;d_z)$ is a constant matrix of linear drag coefficients, while the matrices $\mbf{E}$ and $\mbf{F}$ are constant rotational drag coefficients \cite{Faessler2018}.

\vspace{-0.5em}
\section{Control} \label{Control}
\subsection{Control Objective and Error Definitions}
The control objective is to generate actuator commands that allow the vehicle to follow a reference trajectory. Denote $\mbf{X}_k$ as an element of the matrix Lie group $SE_2(3)$
\begin{align}
	\mbf{X}_k = \bma{ccc} \mbf{C}_{{ab}_k} & \mbf{v}^{zw/a}_{a_k} & \mbf{r}^{zw}_{a_k} \\ \mbf{0} & 1 & 0 \\ \mbf{0} & 0 & 1 \ema \in SE_2(3),
\end{align}
Denote the desired reference frame $\mathcal{F}_r$, the reference attitude, velocity, and position trajectories as $\mbf{C}_{ar}$, $\mbf{v}^{z_rw/a}_a$, and $\mbf{r}^{z_rw}_a$, respectively. The current states and reference states at timestep $k$ are placed into $\mbf{X}_k \in SE_2(3)$ and $\mbf{X}^r_k \in SE_2(3)$, respectively. The tracking error is defined using a left-invariant error \cite{Barrau2017a}
\begin{align}
	\delta\mbf{X}_k = \mbf{X}^{r^{-1}}_k\mbf{X}_k = \bma{ccc} \delta\mbf{C}_k & \delta\mbf{v}_k & \delta\mbf{r}_k \\ \mbf{0} & 1 & 0 \\ \mbf{0} & 0 & 1 \ema \in SE_2(3), \label{eq411}
\end{align}
where $\delta\mbf{C}_k = \mbf{C}_{ar_k}^\trans\mbf{C}_{ab_k}$, $\delta\mbf{v}_k = \mbf{C}_{ar_k}^\trans(\mbf{v}^{zw/a}_{a_k} - \mbf{v}^{z_rw/a}_{a_k})$, and $\delta\mbf{r}_k = \mbf{C}_{ar_k}^\trans(\mbf{r}^{zw}_{a_k} - \mbf{r}^{z_rw}_{a_k})$. The tracking error is expressed in terms of the Lie algebra $\mathfrak{se}_2(3)$ as $\delta \mbf{X}_k = \exp (\delta \mbs{\xi}_k^\wedge)$ where $\delta \mbs{\xi}_k = [\delta \mbs{\xi}_k^{\phi ^\trans} \; \delta \mbs{\xi}_k^{v ^\trans} \; \delta \mbs{\xi}_k^{r ^\trans}]^\trans \in \mathbb{R}^9$ are the matrix Lie algebra representations of attitude, velocity, and position error. See the Appendix for further details about $SE_2(3)$ and $\mathfrak{se}_2(3)$.

Inner-loop control of the rotational dynamics requires a tracking error definition for the angular velocity. The angular momentum tracking error, $\delta\mbf{h}_k \in \mathbb{R}^3$, is defined as
\begin{align}
	\delta\mbf{h}_k = \delta\mbf{C}_k\mbf{h}^{\mathcal{B}z/a}_{b_k} - \mbf{h}^{\mathcal{B}_rz/a}_{r_k} \in \mathbb{R}^3, \label{eq416}
\end{align}
where $\mbf{h}^{\mathcal{B}_rz/a}_{r_k} = \mbf{J}^{\mathcal{B}_rz}_r\mbs{\omega}^{ra}_{r_k}$ is the reference angular momentum and $\mbf{h}^{\mathcal{B}z/a}_{b_k} = \mbf{J}^{\mathcal{B}z}_b\mbs{\omega}^{ba}_{b_k}$ is the true angular momentum. Note that angular momentum is used in the tracking error definition because it results in simpler Jacobians. The control objective of both outer and inner-loop controllers is to drive the tracking error to zero such that $\delta\mbs{\xi}_k = \mbf{0}$ and $\delta\mbf{h}_k = \mbf{0}$.

\vspace{-1em}
\subsection{Overview of Control Structure}
\begin{figure}[t]
	\centering
	\vspace{1em}
	\centerline{\includegraphics[width=\columnwidth]{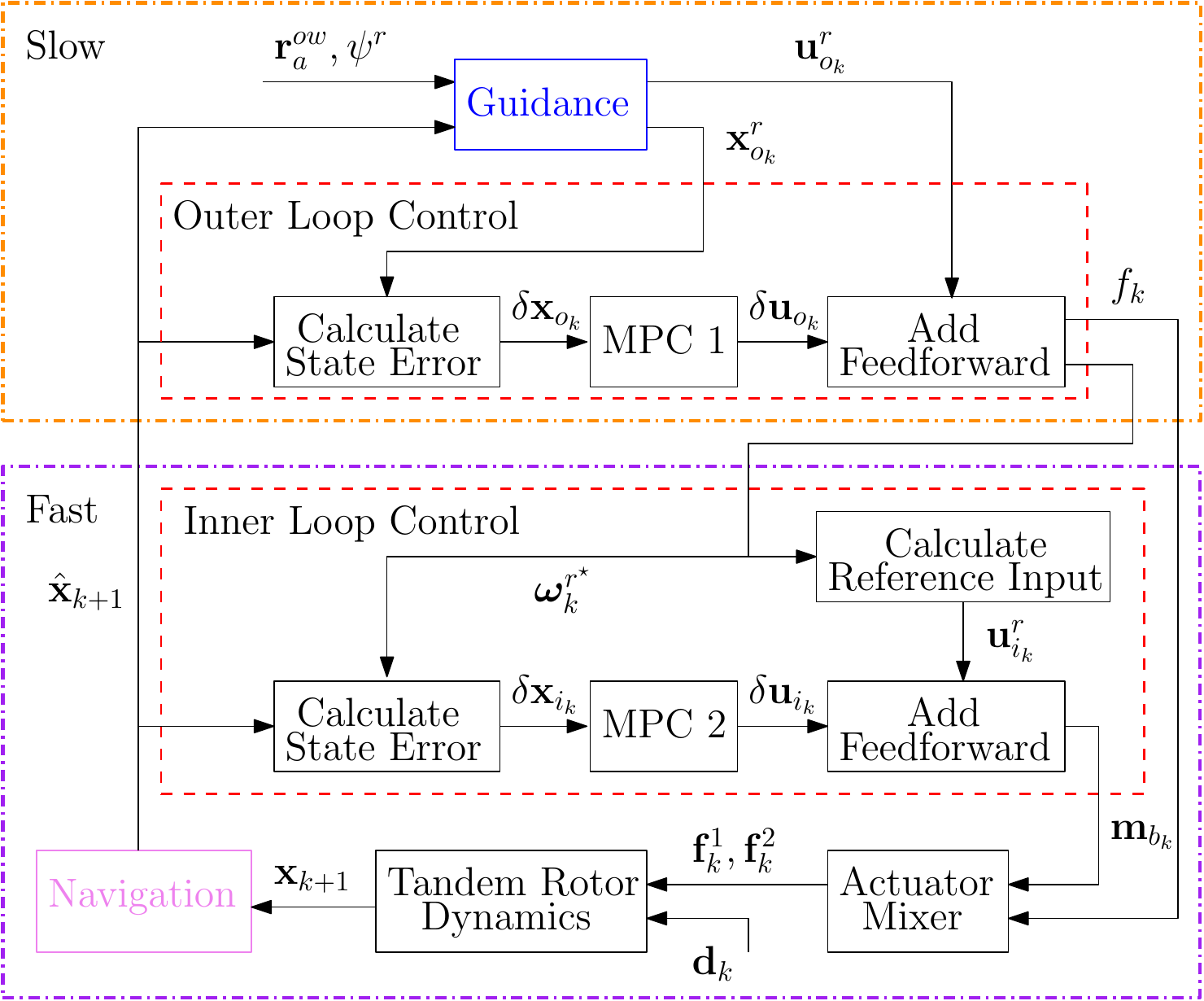}}
	\vspace{-0.5em}
	\caption{Proposed CMPC structure showing multi-timescale approach.}
	\label{fig1a}
	\vspace{-1.5em}
\end{figure}

The proposed cascaded control structure is shown in Fig.~\ref{fig1a}. First, a guidance law generates a reference trajectory in real-time. As demonstrated in \cite{Ahmed2022}, the differentially flat property \cite{Sferrazza2016} of the tandem-rotor system is leveraged such that the full reference state $\mbf{x}^r_o$, and control input $\mbf{u}^r_o$, trajectories for the outer-loop controller can be generated from a desired final position $\mbf{r}^{ow}_a$, and heading $\psi^r$ \cite{Faessler2018}.

Next, the outer-loop state error, $\delta\mbf{x}_{o_k} = \delta\mbs{\xi}_k \in \mathbb{R}^9$, is calculated using $\eqref{eq411}$. The outer-loop \textit{MPC 1} algorithm operates on $\delta\mbf{x}_{o_k}$ to produce the sequence of optimal feedback control inputs, $\delta\mbf{u}^\star_{o_k}$. The first element of $\delta\mbf{u}^\star_{o_k}$ is extracted to obtain $\delta\mbf{u}_{o_k} = [\delta f_k \; \delta\mbs{\omega}_k^\trans]^\trans \in \mathbb{R}^4$, which is added to the feedforward reference control input, $\mbf{u}^r_{o_k} = [f^r_k \; \mbs{\omega}^{ra^\trans}_{r_k}]^\trans~\in~\mathbb{R}^4$, to produce the total thrust force control input, $f_k$, and angular velocity inner-loop input, $\mbs{\omega}^{ba,\text{cmd}}_{b_k}$,
\begin{align}
	\mbf{u}_{o_k} = \bma{c} f_k \\ \mbs{\omega}^{ba,\text{cmd}}_{b_k} \ema = \bma{c} f^{r}_k + \delta f_k \\
	\delta\mbf{C}^\trans_k \mbs{\omega}^{ra}_{r_k} + \delta\mbs{\omega}_{k} \ema \in \mathbb{R}^4. \label{eq420}
\end{align}
Note that \eqref{eq420} can be applied to all elements of $\delta\mbf{u}^\star_{o_k}$ to produce an optimal angular velocity input sequence, $\mbs{\omega}^{r^\star}_k$.

Next, $\mbs{\omega}^{r^\star}_k$ is used as a reference angular velocity trajectory. The differentially flat property of the system is again leveraged to generate a reference rotor torque trajectory, $\mbf{u}^r_{i_k} = \mbf{m}^r_{r_k} \in \mathbb{R}^3$. The inner-loop state error, $\delta\mbf{x}_{i_k} = \delta\mbf{h}_k \in \mathbb{R}^3$, is calculated using \eqref{eq416}. The inner-loop \textit{MPC~2} algorithm operates on $\delta\mbf{x}_{i_k}$ to produce the feedback control input, $\delta\mbf{u}_{i_k} = \delta\mbf{m}_{k} \in \mathbb{R}^3$, which is added to the feedforward reference control input, $\mbf{u}^r_{i_k}$, to produce the total rotor torque control input
\begin{align}
	\mbf{u}_{i_k} = \mbf{m}_{b_k} = \delta\mbf{C}^\trans_k\mbf{m}^r_{r_k} + \delta\mbf{m}_k \in \mathbb{R}^3. \label{eq421}
\end{align}
Because \textit{MPC~2} runs faster than \textit{MPC~1}, the angular velocity reference trajectory $\mbf{u}^{r^\star}_k$ is held by the Inner Loop block until the Outer Loop block produces new output. However, the Inner Loop block steps through $\delta\mbf{u}^{r^\star}$ and re-linearizes the inner-loop dynamics at every inner-loop time step with the current reference angular velocity.

An actuator mixer \cite{Lee2010} is used to map the thrust input from \eqref{eq420} and the torque input from \eqref{eq421} to front and rear rotor force components, $\mbf{f}^1_k$ and $\mbf{f}^2_k$, respectively.

\vspace{-0.5em}
\subsection{Linearization of Dynamics}
Using a Taylor series expansion, first-order approximations for the attitude, velocity, and position errors can be made \cite{Hartley2019}, 
\begin{subequations}\label{eq439}
	\begin{align}
		\delta\mbf{C} &= \exp(\delta\mbs{\xi}^{\phi^\times}) \approx \eye + \delta\mbs{\xi}^{\phi^\times},\\
		\delta\mbf{v} = \mbftilde{J}(\delta\mbs{\xi}^\phi)\delta\mbs{\xi}^v &\approx \delta\mbs{\xi}^v,\quad
		\delta\mbf{r} = \mbftilde{J}(\delta\mbs{\xi}^\phi)\delta\mbs{\xi}^r \approx \delta\mbs{\xi}^r,
	\end{align}
\end{subequations}
where $\delta\mbs{\xi}$ is small, and $\mbftilde{J}(\cdot)$ is the $SO(3)$ left Jacobian defined in the Appendix. Using \eqref{eq439}, and the error definitions \eqref{eq411} and \eqref{eq420}, the continuous-time outer-loop dynamics are linearized about the reference trajectory yielding $\delta\mbfdot{x}_o = \mbf{A}_o\delta\mbf{x}_o + \mbf{B}_o\delta\mbf{u}_o$ where $\delta\mbf{x}_o = [\delta\mbs{\xi}^{\phi^\trans} \;\; \delta\mbs{\xi}^{v^\trans} \;\; \delta\mbs{\xi}^{r^\trans}]^\trans$, $\delta\mbf{u}_o = [\delta f \;\; \delta\mbs{\omega}^\trans]^\trans$, 
\begin{align}
	\mbf{A}_o = \bma{ccc} \mbf{0} & \mbf{0} & \mbf{0} \\ \mbf{A}_{21} & \mbf{A}_{22} & \mbf{0} \\ \mbf{0} & \eye & \mbf{A}_{33} \ema,\quad
	\mbf{B}_o = \bma{cc} \mbf{0} & \eye \\ \mbf{B}_{21} & \mbf{0} \\ \mbf{0} & \mbf{0} \ema, \label{eq440}
\end{align}
and 
\begin{align}
	\begin{split}
		\mbf{A}_{21} &= \f{1}{m_\mathcal{B}}\left(\left(\mbf{D}\mbf{C}_{ar}^\trans\mbf{v}^{z_rw/a}_a\right)^\times \right.\\ &\qquad \left. - \mbf{D}\left(\mbf{C}_{ar}^\trans\mbf{v}^{z_rw/a}_a\right)^\times + \left(f^r \eye_3\right)^\times \right),
	\end{split} \label{eq441}\\
	\mbf{A}_{22} &= -\mbs{\omega}^{ra^\times}_r - \frac{\mbf{D}}{m_\mathcal{B}},\;\;
	\mbf{A}_{33} = -\mbs{\omega}^{ra^\times}_r,\;\;
	\mbf{B}_{21} = -\frac{\eye_3}{m_\mathcal{B}}. \nonumber
\end{align}
Notice that the outer-loop Jacobians from \eqref{eq440} only depend on reference quantities. Importantly, when $\mbf{D} = \mbf{0}$, the Jacobians depend only on control inputs $f^r$ and $\mbs{\omega}^{ra}_r$. A partial derivation of the outer-loop dynamics can be found in the Appendix.

Similarly, the continuous-time inner-loop equations of motion are linearized about the reference angular velocity trajectory using \eqref{eq411}, \eqref{eq416}, \eqref{eq421}, and \eqref{eq439}, yielding 
\begin{align}
	\begin{split}
		\delta\mbfdot{h} &=  \Bigg(\left(\mbf{E}\mbf{C}_{ar}^\trans\mbf{v}^{z_rw/a}_a\right)^\times - \mbf{E}\left(\mbf{C}_{ar}^\trans\mbf{v}^{z_rw/a}_a\right)^\times \\
		&\qquad + \left(\mbf{F}\mbs{\omega}^{ra}_r\right)^\times - \mbf{F}\mbs{\omega}^{ra^\times}_r\Bigg)\delta\mbs{\xi}^\phi - \mbf{E}\delta\mbs{\xi}^v \\
		&\qquad + \left(-\mbs{\omega}^{ra^\times}_r - \mbf{F}\mbf{J}^{\mathcal{B}z_r^{-1}}_r\right)\delta\mbf{h} + \delta\mbf{m}.
	\end{split}
\end{align}
Notice that when the drag terms $\mbf{E} = \mbf{F} = \mbf{0}$, the linearized dynamics can be written as $\delta\mbfdot{x}_i = \mbf{A}_i\delta\mbf{x}_i + \mbf{B}_i\delta\mbf{u}_i$ where $\delta\mbf{x}_i = \delta\mbf{h}$, $\delta\mbf{u}_i = \delta\mbf{m}$, 
\begin{align}
	\mbf{A}_i = -\mbs{\omega}^{ra^\times}_r, \qquad \mbf{B}_i = \eye. \label{eq14}
\end{align}
From \eqref{eq14}, it can be seen that a reference angular velocity trajectory $\mbs{\omega}^{ra}_r$ is required in $\mbf{A}_i$ to linearize the inner-loop dynamics. Recall that the outer-loop MPC generates a sequence of optimal control inputs, $\mbf{u}^\star_k$. The sequence of optimal angular velocity commands, $\mbs{\omega}^\star_k$ are extracted from $\mbf{u}^\star_k$ yielding $\mbs{\omega}^\star_k = [\mbs{\omega}^{ra^\trans}_{r_{0|k}} \; 	\mbs{\omega}^{ra^\trans}_{r_{1|k}} \; \cdots \; \mbs{\omega}^{ra^\trans}_{r_{N-1|k}} ]^\trans$.

Instead of using the reference angular velocity from the guidance, the inner-loop dynamics are linearized about $\mbs{\omega}^\star_k$. The resulting linearization more accurately represents the desired trajectory. From the definition of $\delta\mbf{m}$ given by \eqref{eq421}, a reference torque input trajectory, $\mbf{m}^{r}_{r_k}$ must also be supplied to the inner-loop MPC. This is accomplished using the differential flatness property of the system and the method shown in \cite{Farrell2008}.

The outer and inner-loop continuous-time linearized systems are then discretized \cite{Farrell2008} using their respective controller timesteps yielding
\begin{align}
	\delta\mbf{x}_{k+1} = \mbf{A}_k \delta\mbf{x}_k + \mbf{B}_k \delta\mbf{u}_k. \label{eq422}
\end{align}

\subsection{Finite Horizon MPC for LTV Systems}
For a discretized-linearized system, the state predictions are $\delta\mbf{x}_{i|k} = \bm{\mathcal{A}}_i(k)\delta\mbf{x}_k + \bm{\mathcal{C}}_i(k)\delta\mbs{\mu}_k,\quad i = 0,\ldots,N,$ where the notation $(\cdot)_{i|k}$ denotes a quantity at time $k+i$ predicted at time $k$, and $N$ is the prediction horizon length. The predicted input sequence is $\delta\mbs{\mu}_k~=~[\delta\mbf{u}_{0|k}^\trans \; \cdots \; \delta\mbf{u}_{N-1|k}^\trans]^\trans$ and the time-varying state transition matrices are \cite{Kouvaritakis2016}
\begin{equation}\label{eq423}
	\begin{aligned}
		\mbc{A}_i(k) &= \overset{\curvearrowleft}{\prod^{i-1}_{j=0}} \mbf{A}_{k+j},\\
		\mbc{C}_i(k) &= \Bigg[\begin{matrix}
			\left(\overset{\curvearrowleft}{\prod^{i-1}_{j=1}} \mbf{A}_{k+j}\right)\mbf{B}_k & \left(\overset{\curvearrowleft}{\prod^{i-1}_{j=2}} \mbf{A}_{k+j}\right)\mbf{B}_{k+1}\end{matrix}\Bigg.\\ &\qquad\qquad
		\Bigg.\begin{matrix} \cdots & \mbf{B}_{k+i-1} & \mbf{0} & \cdots & \mbf{0} \end{matrix}\Bigg],
	\end{aligned}
\end{equation}
where $\overset{\curvearrowleft}{\prod}$ indicates that successive terms in the sequence are left-multiplied. The predicted state sequence over the prediction horizon, $\delta\mbs{\chi}_k = [\delta\mbf{x}_{1|k}^\trans \; \cdots \; \delta\mbf{x}_{N|k}^\trans]^\trans$, can be written compactly as \cite{Kouvaritakis2016}
\begin{align}
	\delta\mbs{\chi}_k = \mbc{S}_k\delta\mbs{\mu}_k + \mbc{M}_k\delta\mbf{x}_k, \label{eq424}
\end{align}
where $\mbc{S}_k = [\mbc{C}_1^\trans(k) \; \cdots \; \mbc{C}_N^\trans(k)]^\trans$ and $\mbc{M}_k = [\mbc{A}_1^\trans(k) \; \cdots \; \mbc{A}_N^\trans(k)]^\trans$. Using a standard finite-horizon quadratic cost function \cite{Boyd2004} with state penalty $\mbf{Q} = \mbf{Q}^\trans \geq 0$, control input penalty $\mbf{R} = \mbf{R}^\trans > 0$, and terminal state penalty $\mbf{P} = \mbf{P}^\trans \geq 0$, the constrained optimization problem can be expressed as a QP in terms of $\delta\mbs{\mu}_k$
\begin{equation}\label{eq431}
	\begin{aligned}
		&\underset{\delta\mbs{\mu}_k}{\text{min}}\quad J(\delta\mbs{\mu}_k) = \frac{1}{2}\delta\mbs{\mu}^\trans_k\mbf{H}_k\delta\mbs{\mu}_k + \mbf{F}^\trans_k\delta\mbs{\mu}_k,\\
		&\text{s.t.}\quad \mbf{G}_k\delta\mbs{\mu}_k \leq \mbf{W}_k + \mbf{T}_k\delta\mbf{x}_k,
	\end{aligned}
\end{equation}
where $\mbf{H}_k = \mbc{S}_k^\trans \mbfbar{Q} \mbc{S}_k + \mbfbar{R}$, $\mbf{F}_k = \mbc{S}_k^\trans \mbfbar{Q} \mbc{M}_k$, $\mbfbar{Q}~=~ \text{diag}\left(\mbf{Q},\ldots,\mbf{Q},\mbf{P}\right)$, $\mbfbar{R} = \text{diag}\left(\mbf{R},\ldots,\mbf{R}\right)$, and $\mbf{G}_k$, $\mbf{W}_k$, and $\mbf{T}_k$ are constraint matrices. The solution to \eqref{eq431} gives the optimal control input sequence $\delta\mbs{\mu}^\star_k$.

\subsection{Constraints}
One of the main benefits of the MPC framework is the ability to embed state and input constraints in the optimization problem. The attitude is limited at the outer-loop level using a combination of an attitude keep-in zone and $\ell_1$-norm constraint on the attitude error \cite{Ahmed2022}. The keep-in zone is defined as
\begin{align}
	\mbf{x}_b^\trans \mbf{C}_{ba}^\trans \mbf{y}_a \geq \cos(\alpha) - \epsilon_1, \label{eqkiz}
\end{align}
where $\alpha$ is the keep-in zone angle and $\epsilon_1$ is a slack variable. By setting $\mbf{x}_b = \mbf{y}_a = [0\;0\;1]^\trans$, the roll and pitch angles are simultaneously constrained. The $\ell_1$ attitude error constraint is defined as
\begin{align}
	\norm{\delta\mbs{\xi}^\phi}_1 \leq \gamma + \epsilon_2, \label{eqell1}
\end{align}
where $\gamma$ is the maximum allowable attitude error and $\epsilon_2$ is a slack variable. Control input constraints at both the outer and inner-loop levels are given by
\begin{align}
	\mbf{u}_\text{min} \leq \mbf{u} \leq \mbf{u}_\text{max}.
\end{align}
Note that the state constraints \eqref{eqkiz} and \eqref{eqell1} are implemented as soft constraints to ensure recursive feasibility. A detailed derivation of the constraint linearization and matrices $\mbf{G}_k$, $\mbf{W}_k$ and $\mbf{T}_k$ can be found in \cite{Ahmed2022}.

\subsection{Inner-loop LTI Assumption}
Notice from \eqref{eq14} that the linearized inner-loop dynamics are a linear-time-varying system because the reference angular velocity in $\mbf{A}_i$ is time-varying. Consider the cascaded control structure shown in Fig.~\ref{fig1a}. Assuming that the outer-loop MPC is tuned such that the control inputs generated by the controller vary slowly, the reference angular velocity trajectory used to linearize the inner-loop dynamics can be approximated as a constant value for each iteration of the outer-loop controller, $\mbs{\omega}^\star_k \approx [ \mbs{\omega}^{ra^\trans}_{r_{0|k}} \; \cdots \; \mbs{\omega}^{ra^\trans}_{r_{0|k}} ]^\trans$. Using this assumption, \eqref{eq14} becomes an LTI system, and the inner-loop dynamics only need to be linearized and discretized once per iteration of the outer-loop controller. The $\mbc{S}_k$ and $\mbc{M}_k$ matrices for the inner-loop MPC are simplified to
\begin{align}
	\mbc{S}_{i_k} = \bma{ccc}
	\mbf{B}_i\\
	\vdots & \ddots\\
	\mbf{A}_i^{N-1}\mbf{B}_i & \cdots & \mbf{B}_i
	\ema, \;
	\mbc{M}_{i_k} = \bma{c}
	\mbf{A}_i\\
	\vdots\\
	\mbf{A}_i^N
	\ema.
\end{align}

\section{Simulations} \label{SimResults}
The proposed cascaded control scheme is tested in simulation using a tandem-rotor helicopter model and the equations of motion from \eqref{eq201} and \eqref{eq202}. 	These simulations test the robustness of the proposed control structure to initial conditions, environmental disturbances, and model uncertainty. The physical properties of the vehicle are given in Tab.~\ref{tab1}. Note that $\mbf{r}^{1z}_b$ and $\mbf{r}^{2z}_b$ refer to the locations of the front and rear rotors with respect to the vehicle's center of mass.

\begin{table}[t]
	\vspace{1em}
	\caption{Tandem-rotor Parameters} \label{tab1}
	\vspace{-1em}
	\begin{center}
		\begin{tabular}{|c|c|c|}
			\hline
			Parameter & Value & Units \\
			\hline
			$m_\mathcal{B}$ & 218 & $\text{kg}$ \\
			\hline
			$\mbf{J}^{\mathcal{B}z}_b$ & $\text{diag}(26.8, 97.6, 87.2)$ & $\text{kg}\cdot\text{m}^2$ \\
			\hline
			$\mbf{r}^{1z}_b$ & $[1.045 \; 0 \; -0.514]^\trans$ & $\text{m}$ \\
			\hline
			$\mbf{r}^{2z}_b$ & $[-0.937 \; 0 \; -0.686]^\trans$ & $\text{m}$ \\
			\hline
		\end{tabular}
	\end{center}
	\vspace{-2em}
\end{table}

For all simulations, a final landing approach is simulated using the reference trajectory and generation procedure from \cite{Ahmed2022}. The nominal initial position, velocity, and attitude is $\mbf{r}^{z_0w}_a = [5\;0\;-0.5]^\trans \text{ m}$, $\mbf{v}^{z_0w/a}_a = [5\;0\;-0.5]^\trans \text{ m/s}$, and $\mbf{C}_{ab_0} = \eye$, respectively. The target position, velocity, and heading are $\mbf{r}^{z_fw}_a = \mbf{0}\text{ m}$, $\mbf{v}^{z_fw/a}_a = \mbf{0}\text{ m/s}$, and $\quad \psi^r = 0 \text{ rad}$, respectively.

In practice a navigation loop would provide state estimates, $\mbfhat{x}_k$, to the guidance and control algorithms. For these simulations, it is assumed the state estimates are accurate such that $\mbf{x}_k = \mbfhat{x}_k$.

\subsection{Conventional Single MPC Structure}
The proposed CMPC structure is compared to the SMPC structure presented in \cite{Ahmed2022}. In the conventional approach, the $SE_2(3)$ error state is augmented with the angular momentum to form a single state-space system $\delta\mbfdot{x}_c = \mbf{A}_c\delta\mbf{x}_c + \mbf{B}_c\delta\mbf{u}_c$ where $\delta\mbf{x}_c = [\delta\mbs{\xi}^{\phi^\trans} \; \delta\mbs{\xi}^{v^\trans} \; \delta\mbs{\xi}^{r^\trans} \; \delta\mbf{h}^\trans]^\trans \in \mathbb{R}^{12}$, $\delta\mbf{u}_c = [\delta f \; \delta\mbf{m}^\trans]^\trans \in \mathbb{R}^4$, 
\begin{align}
	\mbf{A}_c = \bma{ccccc} \mbf{0} & \mbf{0} & \mbf{0} & \mbf{A}_{14} \\ \mbf{A}_{21} & \mbf{A}_{22} & \mbf{0} & \mbf{0} \\ \mbf{0} & \eye & \mbf{A}_{33} & \mbf{0} \\ \mbf{A}_{41} & \mbf{A}_{42} & \mbf{0} & \mbf{A}_{44} \ema,\;
	\mbf{B}_c = \bma{cc} \mbf{0} & \mbf{0} \\ \mbf{B}_{21} & \mbf{0} \\ \mbf{0} & \mbf{0} \\ \mbf{0} & \eye \ema, \label{eq450}
\end{align}
and 
\begin{subequations}
	\begin{align}
		\mbf{A}_{14} &= \mbf{J}^{\mathcal{B}z^{-1}}_b, \qquad \mbf{A}_{44} = -\mbs{\omega}^{ra^\times}_r -\mbf{F}\mbf{J}^{\mathcal{B}z_r^{-1}}_r,\\
		\begin{split}
			\mbf{A}_{41} &= \left(\mbf{E}\mbf{C}_{ar}^\trans\mbf{v}^{z_rw/a}_a\right)^\times - \mbf{E}\left(\mbf{C}_{ar}^\trans\mbf{v}^{z_rw/a}_a\right)^\times \\
			&\qquad + \left(\mbf{F}\mbs{\omega}^{ra}_r\right)^\times - \mbf{F}\mbs{\omega}^{ra^\times}_r.
		\end{split}
	\end{align}
\end{subequations}
The remaining terms in $\mbf{A}_c$ and $\mbf{B}_c$ are as defined in \eqref{eq441}.

\begin{table}[t]
	\vspace{1em}
	\caption{MPC Parameters Used in Simulations} \label{tab2}
	\vspace{-1em}
	\begin{center}
		\begin{tabular}{|c|c|c|c|}
			\hline
			Parameter & Outer CMPC & Inner CMPC & SMPC \\
			\hline
			$\Delta t \text{ (s)}$ & 0.10 & 0.02 & 0.02 \\
			\hline
			$N$ & $48$ & $10$ & $48$ \\
			\hline
			$N_u$ & $10$ & $5$ & $10$ \\
			\hline
			$N_c$ & $10$ & $5$ & $10$ \\
			\hline
			$T_\text{total} \text{ (s)}$ & $28.8$ & $0.2$ & $5.76$ \\
			\hline
			\multirow{2}*{$\mbf{Q}$} & $10\cdot\text{diag}(100\cdot\eye,$ & \multirow{2}*{$1000\cdot\eye$} & $10\cdot\text{diag}(100\cdot\eye, \eye,$ \\
			& $\quad \eye, 10\cdot\eye)$ & & $\quad 10\cdot\eye, \eye)$ \\
			\hline
			$\mbf{P}$ & $1\cdot\mbf{Q}$ & $1\cdot\mbf{Q}$ & $1\cdot\mbf{Q}$ \\
			\hline
			$\mbf{R}$ & $\text{diag}(0.001,1,1,1)$ & $\eye$ & $\text{diag}(0.001,1,1,1)$ \\
			\hline
		\end{tabular}
	\end{center}
	\vspace{-2em}
\end{table}

\subsection{MPC Parameters}
The base controller timestep, $\Delta t$, prediction horizon length, $N$, control horizon length, $N_u$, constraint horizon length, $N_c$, and total prediction time, $T_\text{total}$, for the CMPC and SMPC structures are shown in Tab.~\ref{tab2}. A non-uniform prediction horizon timestep \cite{Ahmed2022, Tan2016} is applied to both the outer-loop CMPC and the SMPC to extend the total prediction time and improve tracking performance without increasing the optimization problem size. 

The keep-in zone is set to $\alpha = 0.14 \text{ rad}$, while the $\ell_1$ attitude error constraint is set to $\gamma = 0.1 \text{ rad}$. The control input constraints are $0 \leq f \leq 3000 \text{ N}$, $-2 \leq \omega^{ba}_{b_i} \leq 2 \text{ rad}\cdot\text{s}^{-1}$, and $-200 \leq m_{b_i} \leq 200 \text{ N}\cdot\text{m}$. For the SMPC loop, the angular velocity constraint is applied as a state constraint on the angular momentum.

The reference trajectory is replanned in real-time during the simulation if the $\ell_1$ constraint is active for longer than $0.4$~s to ensure the linearization of the dynamics remain valid. This value is tuned empirically to ensure the trajectory is not continuously replanned.

\begin{figure*}[t]
	\vspace{0.5em}
	\centering
	\includegraphics[width=\textwidth]{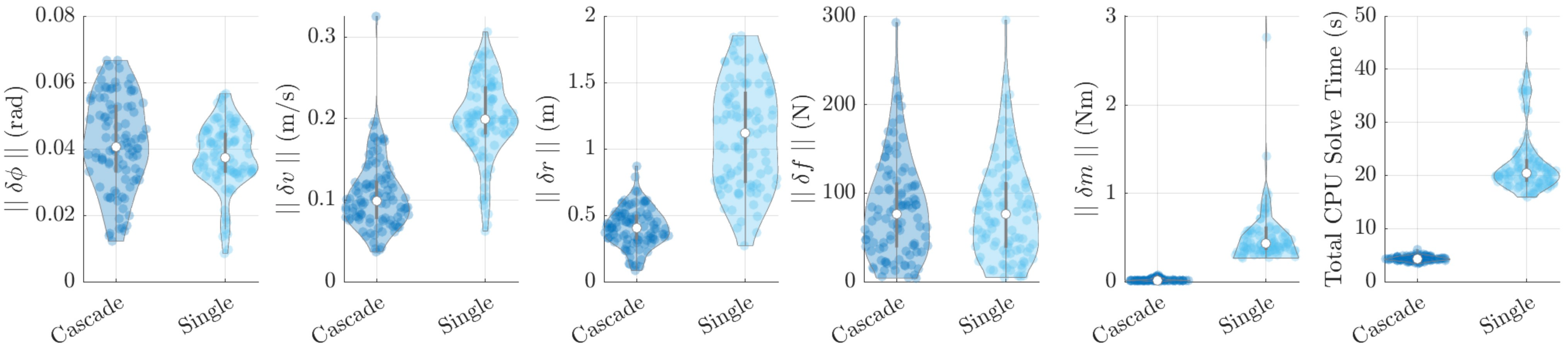}
	\vspace{-1.5em}
	\caption{State/input RMSE values, and total CPU solve time across $100$ Monte-Carlo simulations comparing the CMPC and SMPC structures.}
	\label{fig7}
	\vspace{-1em}
\end{figure*}

\begin{figure}[t]
	\centering
	\centerline{\includegraphics[width=\columnwidth]{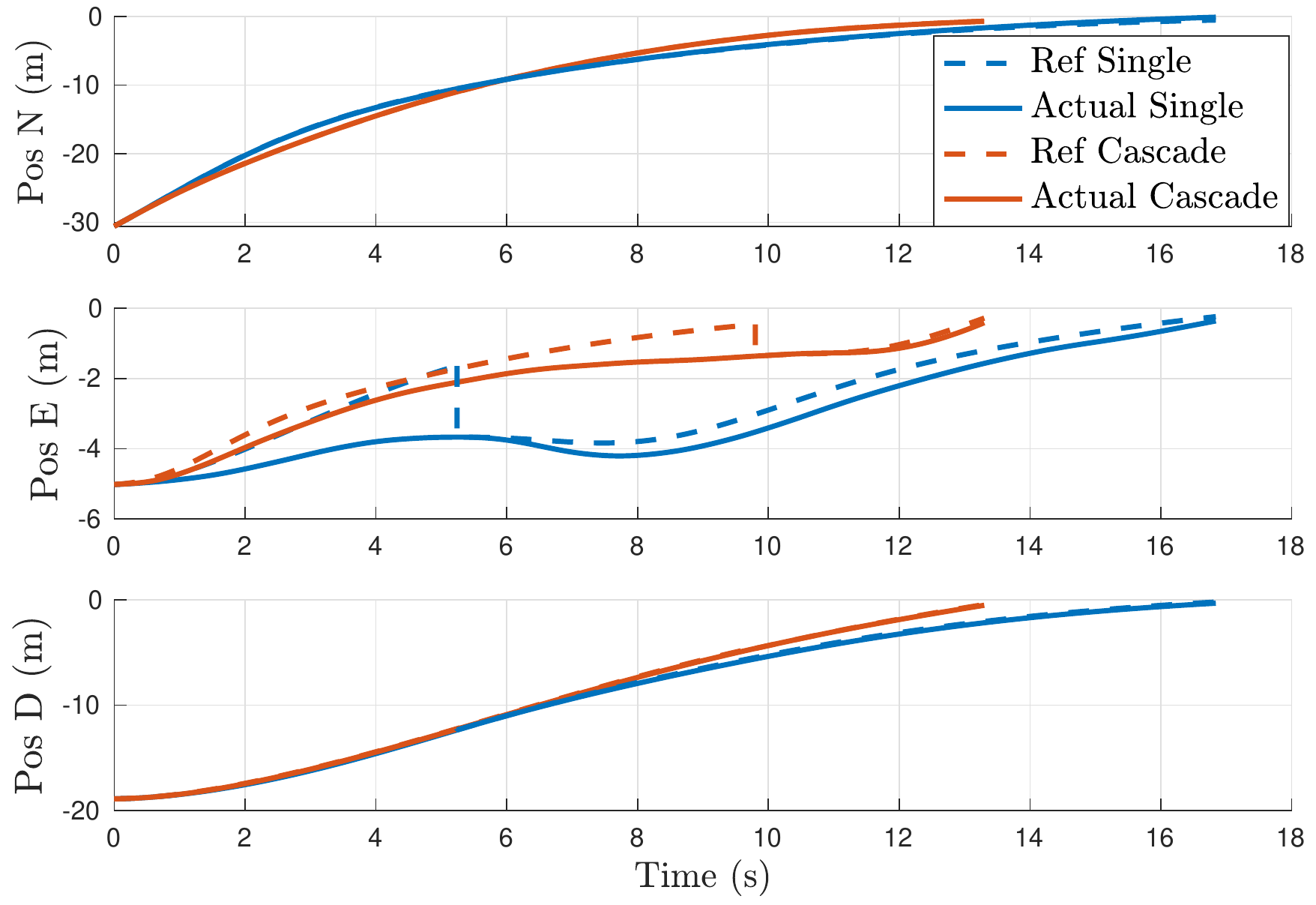}}
	\vspace{-1em}
	\caption{Position trajectory for a single simulation comparing the CMPC and SMPC structures.}
	\label{fig8}
	\vspace{-1.5em}
\end{figure}

\subsection{Monte-Carlo Parameters}
The initial state of the helicopter is randomized such that $\mbs{\phi}_0 = \mbf{0} + \mbf{w}_1 \text{ rad}$, $\mbf{r}^{z_0w}_a = [-30 \; -5 \; -20]^\trans + \mbf{w}_2 \text { m}$, $\mbf{v}^{z_0w/a}_a = [5 \; 0 \; 0.5]^\trans + \mbf{w}_3 \text{ m/s}$, and $\mbs{\omega}^{b_0a}_b = \mbf{0} + \mbf{w}_4 \text{ rad/s}$, where $\mbs{\phi}_0 = \log_{SO(3)}(\mbf{C}_{ab_0})$, $\mbf{w}_1 \sim \mathcal{N}(0,0.116^2\eye)$,  $\mbf{w}_2 \sim \mathcal{N}(0,\eye)$,  $\mbf{w}_3 \sim \mathcal{N}(0,0.333^2\eye)$,  and $\mbf{w}_4 \sim \mathcal{N}(0,0.029^2\eye)$. 

The nominal wind speed is set to $\mbf{v}^{sw/a}_a = [0 \; -5 \;\; 0]^\trans + \mbf{w}_5 \text{ m/s}$ where $\mbf{w}_5 \sim \mathcal{N}(0,1.667^2\eye)$. Wind gusts, $\mbf{v}^{gw/a}_a$, are generated using the Dryden model \cite{Moorhouse1982} and randomized such that $W_0 \sim \mathcal{N}(10 \text{ m/s},1)$, where $W_0$ is the low altitude intensity of the Dryden model.

The estimated mass, $\hat{m}_\mathcal{B}$, and inertia matrix, $\mbfhat{J}^{\mathcal{B}z}_b$, are perturbed from their true values such that $\hat{m}_\mathcal{B} = m_\mathcal{B} + w_5$, and $\mbfhat{J}^{\mathcal{B}z}_{\hat{b}} = \mbf{C}^\trans_{b\hat{b}} \mbf{J}^{\mathcal{B}z}_b \mbf{C}_{b\hat{b}}$, 
where $w_5 \sim \mathcal{N}(0,10^2)$, and $\mbf{C}_{b\hat{b}} = \exp_{SO(3)}(\mbs{\phi})$ is a perturbation DCM, where $\mbs{\phi} \sim \mathcal{N}(0,0.044^2\eye)$.

\vspace{-1em}
\subsection{Results}
The distribution of root-mean-square error (RMSE) results for $100$ Monte-Carlo runs  is given by Fig.~\ref{fig7} in the form of violin plots. In each case, the helicopter is able to reach the target position. The attitude tracking performance of the cascaded controller is similar to that of the single controller. However, the cascaded controller shows dramatically improved velocity and position tracking by $50$\% and $64$\% respectively. Of particular note is the cascaded controller's tighter distribution of position tracking error, indicating that it is better able to reject the various disturbances. 

Both controllers performed nearly identically in terms of thrust input tracking, however the cascaded controller was superior in terms of torque input tracking. Not only is the cascaded controller able to reduce the torque input tracking error by $96$\%, it is noticeably more consistent.

The position trajectory for a single Monte-Carlo simulation is shown in Fig.~\ref{fig8}. The dotted and solid lines indicate the reference and realized trajectories, respectively. The step changes in the reference trajectory occur when it is replanned. A large deviation in East-direction tracking is observed due to a large wind disturbance. The cascaded algorithm is less affected by the disturbance and is able to better track the reference trajectory. Additionally, the cascaded algorithm reaches the target $21$\% faster on average.

The improved performance of the CMPC architecture is attributable to the ability to independently tune the outer and inner-loop MPC algorithms. The dynamics of the tandem-rotor helicopter are multi-timescale. The translational dynamics operate on a much larger time constant than the rotational dynamics. Although a non-uniformly spaced prediction horizon is used in the SMPC controller to account for the different timescales, the performance is still compromised. In the cascaded approach, the outer and inner-loop MPC penalty matrices are tuned to better suit the translational and rotational dynamics separately.

An additional benefit of the cascaded approach is the drastic reduction in computational effort. The total CPU solve time during each Monte-Carlo simulation is shown in the right-most plot of Fig.~\ref{fig7}. Note that the cascaded approach requires solving two QPs at each iteration, but that the problem size for each is reduced since $\delta\mbf{x}_o \in \mathbb{R}^9$ and $\delta\mbf{x}_i \in \mathbb{R}^3$ compared to $\delta\mbf{x}_c \in \mathbb{R}^{12}$. As a result, the total CPU time of the cascaded controller is reduced by $79$\%. More importantly, the cascaded controller solve time is far more consistent with a standard deviation of $0.42$ s compared to $6.06$ s. No task execution overruns are observed in the outer and inner-loops with the cascaded control structure, which provides evidence that the proposed control structure is suitable for real-time implementation.

\vspace{-1em}
\section{Conclusion}
A cascaded MPC approach for a tandem-rotor helicopter is presented in this paper. In the outer-loop, an $SE_2(3)$ error definition is used to linearize the process model about a reference trajectory. In the inner-loop, the rotational dynamics are linearized using the optimal angular velocity sequence from the outer-loop MPC. Monte-Carlo simulations demonstrate robustness to initial conditions, model uncertainty, and environmental disturbances while imposing attitude and control input constraints. Additionally, the cascaded MPC structure is shown to outperform the tracking performance of a single MPC structure while using significantly less computational resources. Future work will focus on hardware implementation and real-world testing.

\section{Appendix}
\subsection{The Group of Double Direct Isometries $SE_2(3)$}
The group of double direct isometries, $SE_2(3)$, introduced in \cite{Barrau2015}, is given by
\begin{align*}
	SE_2(3) = \mbf{X} = \bma{ccc} \mbf{C} & \mbf{v} & \mbf{r} \\ \mbf{0} & 1 & 0 \\ \mbf{0} & 0 & 1 \ema \in \mathbb{R}^{5\times5},
\end{align*}
where $\mbf{C} \in SO(3)$ is a DCM, $\mbf{v} \in \mathbb{R}^3$ is the velocity, and $\mbf{r} \in \mathbb{R}^3$ is the position. The associated Lie algebra, $\mathfrak{se}_2(3)$, is given by
\begin{align*}
	\mathfrak{se}_2(3) = \mbs{\Xi} = \mbs{\xi}^\wedge \in \mathbb{R}^{5\times5},
\end{align*}
where $\mbs{\xi} \in \mathbb{R}^9$, and
\begin{align*}
	\mbs{\xi}^\wedge = \bma{c} \mbs{\xi}^\phi \\ \mbs{\xi}^v \\ \mbs{\xi}^r \ema^\wedge = \bma{ccc} \mbs{\xi}^{\phi^\times} & \mbs{\xi}^v & \mbs{\xi}^r \\ \mbf{0} & 0 & 0 \\ \mbf{0} & 0 & 0 \ema.
\end{align*}
The exponential map from $\mathfrak{se}_2(3)$ to $SE_2(3)$ is
\begin{align*}
	\exp\left(\mbs{\xi}^\wedge\right) = \bma{ccc} \exp\left(\mbs{\xi}^{\phi^\times}\right) & \mbftilde{J}\left(\mbs{\xi}^\phi\right)\mbs{\xi}^v & \mbftilde{J}\left(\mbs{\xi}^\phi\right)\mbs{\xi}^r \\ \mbf{0} & 0 & 0 \\ \mbf{0} & 0 & 0 \ema,
\end{align*}
while the logarithmic map from $SE_2(3)$ to $\mathfrak{se}_2(3)$ is
\begin{align*}
	\log(\mbf{X}) = \bma{ccc} \log_{SO(3)}(\mbf{C}) & \mbftilde{J}\left(\mbs{\xi}^\phi\right)\inv\mbf{v} & \mbftilde{J}\left(\mbs{\xi}^\phi\right)\inv\mbf{r} \\ \mbf{0} & 0 & 0 \\ \mbf{0} & 0 & 0 \ema,
\end{align*}
where the left Jacobian and its inverse are given by \cite{Hartley2019}
\begin{align*}
	\mbftilde{J}(\mbs{\phi}) &= \f{\sin\phi}{\phi}\eye + \left(1 - \f{\sin\phi}{\phi}\right)\mbf{aa}^\trans + \f{1 - \cos\phi}{\phi}\mbf{a}^\times, \\
	\mbftilde{J}(\mbs{\phi})\inv &= \f{\phi}{2}\cot\f{\phi}{2}\eye + \left(1 - \f{\phi}{2}\cot\f{\phi}{2}\right)\mbf{aa}^\trans - \f{\phi}{2}\mbf{a}^\times.
\end{align*}

\subsection{Attitude Linearization}
A partial derivation of the tandem-rotor dynamics is given here. The time rate of change of the attitude error is found by differentiating \eqref{eq201} as
\begin{align*}
	\delta\mbfdot{C} &= \mbfdot{C}_{ar}^\trans \mbf{C}_{ab} + \mbf{C}_{ar}^\trans \mbfdot{C}_{ab} \nonumber \\
	&= -\mbs{\omega}^{ra^\times}_r \delta\mbf{C} + \delta\mbf{C} \mbs{\omega}^{ba^\times}_b. \label{eq4231_1}
\end{align*}
Recall from \eqref{eq420} that the angular velocity can be written as  $\mbs{\omega}^{ba}_b = \delta\mbf{C}^\trans\mbs{\omega}^{ra}_r + \delta\mbs{\omega}$. Therefore,
\begin{align*}
	\delta\mbfdot{C} &= -\mbs{\omega}^{ra^\times}_r \delta\mbf{C} + \delta\mbf{C} \left(\delta\mbf{C}^\trans\mbs{\omega}^{ra}_r + \delta\mbs{\omega}\right)^\times. 
\end{align*}
Linearizing using \eqref{eq439} and expanding,
\begin{align*}
	\begin{split}
		\delta\mbsdot{\xi}^{\phi^\times} &= -\mbs{\omega}^{ra^\times}_r \left(\eye + \delta\mbs{\xi}^{\phi^\times}\right) \\
		&\qquad + \left(\eye + \delta\mbs{\xi}^{\phi^\times}\right)\Bigg(\left(\eye - \delta\mbs{\xi}^{\phi^\times}\right)\mbs{\omega}^{ra}_r + \delta\mbs{\omega}\Bigg)^\times
	\end{split} \\
	\begin{split}
		&= -\mbs{\omega}^{ra^\times}_r \delta\mbs{\xi}^{\phi^\times} - \left(\delta\mbs{\xi}^{\phi^\times}\mbs{\omega}^{ra}_r\right)^\times + \delta\mbs{\omega}^\times \\
		&\qquad + \delta\mbs{\xi}^{\phi^\times}\mbs{\omega}^{ra^\times}_r - \delta\mbs{\xi}^{\phi^\times}\left(\delta\mbs{\xi}^{\phi^\times}\mbs{\omega}^{ra}_r\right)^\times + \delta\mbs{\xi}^{\phi^\times}\delta\mbs{\omega}^\times.
	\end{split}
\end{align*}
Dropping higher order terms,
\begin{align*}
	\begin{split}
		\delta\mbsdot{\xi}^{\phi^\times} &= \delta\mbs{\omega}^\times + \delta\mbs{\xi}^{\phi^\times}\mbs{\omega}^{ra^\times}_r -\mbs{\omega}^{ra^\times}_r \delta\mbs{\xi}^{\phi^\times} - \left(\delta\mbs{\xi}^{\phi^\times}\mbs{\omega}^{ra}_r\right)^\times.
	\end{split}
\end{align*}
Making use of the identity $\mbf{u}^\times\mbf{v}^\times - \mbf{v}^\times\mbf{u}^\times = \left(\mbf{u}^\times\mbf{v}\right)^\times$ where $\mbf{u}, \mbf{v} \in \mathbb{R}^3$,
\begin{align*}
	\delta\mbsdot{\xi}^{\phi^\times} &= \delta\mbs{\omega}^\times + \left(\delta\mbs{\xi}^{\phi^\times}\mbs{\omega}^{ra}_r\right)^\times - \left(\delta\mbs{\xi}^{\phi^\times}\mbs{\omega}^{ra}_r\right)^\times = \delta\mbs{\omega}^\times.
\end{align*}
Uncrossing both sides, $\delta\mbsdot{\xi}^\phi = \delta\mbs{\omega}$. A similar linearization process is followed for the velocity, position, and rotational dynamics to get $\delta\mbsdot{\xi}^v$, $\delta\mbsdot{\xi}^r$, and $\delta\mbsdot{\omega}$.

\bibliographystyle{IEEEtran}
\bibliography{IEEEabrv,lcss.bib}

\end{document}